\documentstyle[12pt,twoside]{amsart}
\title{The canonical modifications by weighted blow-ups}
\author{Shihoko Ishii}

\newcommand{\bZ}{{\Bbb Z}}
\newcommand{\bR}{{\Bbb R}}
\newcommand{\bN}{{\Bbb N}}
\newcommand{\bC}{{\Bbb C}}
\newcommand{\bQ}{{\Bbb Q}}

\newcommand{\D}{{\Delta}}

\newcommand{\bp}{{\bold p}}
\newcommand{\bq}{{\bold q}}
\newcommand{\ba}{{\bold a}}
\newcommand{\bb}{{\bold b}}
\newcommand{\be}{{\bold e}}
\newcommand{\bt}{{\bold t}}
\newcommand{\bl}{{\bold 1}}
\newcommand{\br}{{\bold r}}
\newcommand{\bw}{{\bold w}}
\newcommand{\GGf}{{\Gamma_+(f)}}
\newcommand{\Gf}{{\Gamma(f)}}
\newcommand{\GGg}{{\Gamma_+(g)}}
\newcommand{\Gg}{{\Gamma(g)}}
\newcommand{\Cf}{{C_{\bl}(f)}}
\newcommand{\Cg}{{C_{\bl}(g)}}
\newcommand{\Cn}{{{\Bbb C}^{n+1}}}
\newcommand{\Dq}{{\Delta ({\bold q})}}

\newtheorem{thm}{Theorem}[section]

\newtheorem{lem}[thm]{Lemma}
\newtheorem{cor}[thm]{Corollary}
\newtheorem{prop}[thm]{Proposition}
\newtheorem{conj}[thm]{Conjecture}

\theoremstyle{definition}
\newtheorem{defn}[thm]{Definition}

\newtheorem{say}[thm]{}
\newtheorem{exmp}[thm]{Example}

\newtheorem{rem}[thm]{Remark}

\theoremstyle{remark}

\begin{document}
\bibliographystyle{amsplain}
\maketitle
\vskip 2truecm

{\bf Abstract.}
In this paper
  we give a criterion
  for an isolated, hypersurface singularity of dimension  $n\ (\geq 2)$
  to have the canonical modification by means of a suitable
  weighted blow-up.
 Then we give a counter example to the following conjecture by
  Reid-Watanabe:
  For a 3-dimensional, isolated, non-canonical, log-canonical
  singularity  $(X,x)$
  of embedded dimension   4,
  there exists an embedding  $(X,x)\subset ({\bC}^4, 0)$  and a
  weight  ${\bw}=(w_0,w_1,\ldots ,w_n)$,
  such that the
  ${\bw}$-blow-up gives the canonical modification of $(X,x)$.

\setcounter{section}{-1}
\section{Introduction}
\label{introduction}

\begin{say}
\label{def of canonical modification}
Throughout this paper
all varieties are defined over the complex number field  $\bC$.
The canonical modification of a singularity  $(X,x)$  is a partial
resolution  $\varphi : Y \to X$  such that  $Y$  admits at worst
canonical singularities, and the canonical divisor  $K_Y$  is
relatively ample with respect to  $\varphi$.
If a canonical modification exists,
then it is unique up to isomorphisms over  $X$.
It is well known that it exists if the minimal model conjecture holds.
For a 2-dimensional singularity  $(X,x)$,
the canonical modification is the RDP-resolution (\cite{L}).
For a 3-dimensional singularity  $(X,x)$,
the canonical modification exists by the affirmative answer
(\cite{M}) to
the minimal model conjecture.
For higher dimensional singularities the existence is not generally
 proved.

  The motivation for writing this paper is the following conjecture by
Miles Reid
  and Kimio Watanabe:
\end{say}

\begin{conj}
\label{R-W-conjecture}
  For a 3-dimensional, isolated, non-canonical, log-canonical
  singularity  $(X,x)$
  of embedded dimension   4,
  there exists an embedding  $(X,x)\subset ({\bC}^4, 0)$  and a
  weight  ${\bw}=(w_0,w_1,\ldots ,w_n)$  such that the
  ${\bw}$-blow-up gives the canonical modification of $(X,x)$.
\end{conj}
  Primarily,
  both Reid and Watanabe brought up stronger versions in different ways:
  Reid required the statement for elliptic singularities,
  not only
  for log-canonical singularities (Conjecture p.306, \cite{Reid});
   Kimio Watanabe required it for all non-canonical,
  log-canonical singularities defined by a non-degenerate
  polynomial without replacing embedding to  $({\bC}^4, 0)$,
  and he also required the weight $\bw$ should be in 95-weights
  listed by Yonemura \cite{Y}  and Fletcher \cite{F},  which
  give the weights of quasi-homogeneous simple K3-singularities

  Tomari \cite{T} showed an affirmative answer for log-canonical
  singuarities of special type.
  Watanabe calculated many examples, and made a list of standard equations
  of log-canonical singularities which admit the canonical modifications
  by weighted blow-up with each weight of the 95's \cite{W2}.

\begin{say}
\label{purpose}

  In this paper
  we say that a weight is the canonical weight,
  if it gives the  weighted blow-up
  which is the canonical modification.
  We give a criterion
  for an isolated, hypersurface singularity of dimension  $n\ (\geq 2)$
  to have the canonical weight in \S 2.
  As a consequence,
  for a non-canonical, log-canonical singularity
  $(X,x) \subset ({\bC}^{n+1},0)$
  defined by a non-degenerate polynomial  $f$
  (definition cf. \cite{Kouchinirenco}),
  a primitive vector  $\bw =(w_0,w_1,\ldots ,w_n)$
  is the canonical weight
  if and only if  $\bw$  is absolutely minimal (i.e. each coordinate
  $w_i$ is the minimal integer) in the essential cone in the dual space
  (Corollary \ref{log-canonical}).
  We can see many singularities for which such vectors actually exist
  (Corollary \ref{surface}$\sim $ Example \ref{typeT}).
  But we also observe in \S 4
  an example for which such a vector does not exist
  and it turns out to be a counter  example opposing
  the Reid-Watanabe's conjecture.
   In the other sections,
  we prepare the formula for coefficients of divisors (in \S 1 )
  and study deformations of isolated singularities (in \S 3).
\end{say}

{\bf Acknowledgments.}
  The author would like to thank Professors Miles Reid and Kimio
  Watanabe for introducing the conjecture,
  stimulating discussions and providing many examples.
  She also expresses her gratitude to Professor
  Masataka Tomari
  whose example helped to find an error of the first draft of
  this paper.
  She is also grateful to Professor Kei-ichi Watanabe
  as well as other members of Waseda Tuesday Seminar,
  for constructive comments and advice in preparation of this paper.


\section{Divisors on toric varieties}
\begin{say}
\label{coefficient of divisor}
  Let  $M$  be the free abelian group  ${\bZ}^r$ $(r\geq 1)$
  and  $N$  be the dual $Hom_{\bZ}(M, {\bZ})$.
  We denote  $M\otimes _{\bZ}{\bR}$  and $N\otimes_{\bZ}{\bR}$  by
  ${M_{\bR}}$  and  $N_{\bR}$, respectively.
  Then $N_{\bR}=Hom_{\bR}(M_{\bR},{\bR})$.
  For a finite fan ${\D}$ in ${N_{\bR}}$,
  we construct the toric variety  $V=T_N({\D})$.
  Denote by  ${\D}(1)$
  the set of primitive vectors ${\bq}=(q_1,\ldots , q_r)\in {\bN}$
  whose rays  ${\bR}_{\geq 0}{\bq}$  belong to ${\D}$ as one-dimensional
  cones.
  For ${\bq}\in {\D}(1)$,
  denote by $D_{\bq}$ the
  corresponding divisor which is denoted by $\overline{orb\
  {\bR}_{\geq 0}{\bq}}$  in  \cite{Oda75}.
  Denote by $U_{\bq}$  the invariant affine open subset which contains
  ${orb\
  {\bR}_{\geq 0}{\bq}}$ as the unique closed orbit.
   Then  $U_{\bq}= {\rom Spec} {\Bbb C}[ {\bq}^{\vee}\cap
  M]$,
 $U_{\bq}\cap D_{\bq}={orb\ {\bR}_{\geq 0}{\bq}}
  ={\rom Spec}  {\Bbb C}[ {\bq}^{\perp}\cap
  M]$,
  and $U_{\bq}\cap D_{\bq}$ is defined in $U_{\bq}$
  by the ideal ${\frak p}_{\bq}$
  which is generated by the elements ${\bold e }\in M$ with ${\bq}
  ({\bold e})>0$.

  Express $\bC[M]$  by  $\bC[x^{\ba}]_{\ba\in M}$,
  where $x^{\ba}=x_1^{a_1}x_2^{a_2}\cdots x_r^{a_r}$
  for ${\bold a}=(a_1,\ldots ,a_r)\in M$.
  For convenience sake,
  we write    $x^{\bold a}\in f$,
  if
 $f=\sum_{{\bq}\in M}{\alpha}_{\bq}x^{\bq},$
   and ${\alpha}_{\bold a}\neq 0$.

\end{say}

\begin{defn}
\label{weight of divisor}
  Under the notation above,
  take ${\bq}\in |\D|$.

  (1) For a regular function  $f$  on  $T_N(\D)$,
   define :
   $${\bq}(f):= \min\{{\bq}({\bold a})|{\bold a}\in M,\
   x^{\bold a}\in f\}$$.

  (2)  For a $\bQ$-divisor
  $D $  on $ T_N(\D)$
   such that
  $mD$  is defined by a regular  function  $f$  on $ T_N(\D)$,
   define:
   $${\bq}(D):= \dfrac{1}{m}{\bq}(f).$$

\end{defn}
\begin{prop}
\label{valuation of function}
  Let  $D\subset V $ be a $\bQ$-principal divisor (i.e. $mD$  is defined by
  a regular function on  $V$).
  Then  $D$  is of the form:
  $$D=D'+\sum_{{\bq}\in \D (1)} {\bq}(D)D_{\bq},$$
  where  $D'$  is an effective divisor
which does not contain any  $D_{\bq}$.
\end{prop}

\begin{pf}
 By the definition of  $\bq(D)$, we may assume that $D$  is defined
  by a regular function $f$  on  $V$.
  Write  $D=D'+ \sum_{{\bq}\in \D (1)} m_{\bq}D_{\bq}$
  such that  $D'$ does not contain any $D_{\bq}$
  and  $m_{\bq}\geq 0$.
  Let  $\zeta$ be the defining function of
  $D_{\bq}$ on  $U_{\bq}$,
  then  $f={\zeta}^{m_{\bq}}h$,
  where  $h\notin (\zeta)$.
  Since  ${\bq}(f)=m_{\bq}{\bq}({\zeta})+{\bq}(h)$  and $\bq(h)=0$,
  it is sufficient to prove that  ${\bq}(\zeta)=1$.
  And this is clear,
  because there is a vector ${\bold e}\in M$
  such that  ${\bq}({\bold e})=1 $
  and  ${\zeta}$  is a generator of  the ideal $\{g\in \bC[{\bq}^{\vee}\cap M]|
  \bq(g)>0\}$.
\end{pf}
\begin{say}
\label{simplicial}
  Now we consider the case that  $\D$ consists of all faces of a
  simplicial cone  $\sigma$  in $N_{\bR}$.
  Let  ${\ba}_1,\ldots , {\ba}_r$  be the primitive vector of the
  one dimensional faces of  $\sigma \cap N$
  and ${\ba}^*_1,\ldots , {\ba}^*_r$
 be the dual system of $\{{\ba}_i\}$'s
(i.e. ${\ba}_j^*\in M_{\bR}$  and  ${\ba}_i({\ba}_j^*)
  =\delta _{ij}$).
  Denote by  $\overline{N}$ the subgroup of  $N$  generated by
  $\{{\ba}_i\}$.
  Then the morphism  $\pi : T_{\overline N}(\D) \to T_N(\D)=V$
  induced by  $({\overline N},\D) \to (N, \D)$  is the quotient
  morphism of  ${\bC}^r$  by the finite group  $N/{\overline N}$,
  and each  ${\ba}_j^*$  belongs to
  ${\overline M}=Hom_{\bZ}({\overline N}, {\bZ})$.
  Denote by $D_{{\ba}_i}$ and  ${\overline D_{{\ba}_i}}$
 the divisors on $T_N(\D)$ and
  $T _{\overline N}(\D)={\bC}^r$ respectively
  which are
corresponding to  ${\ba}_i$.
  Let  $r_i$  be the ramification index of  $\pi$  at
  ${\overline D_{{\ba}_i}}$.
\end{say}
\begin{lem}
\label{dual}
 Let ${\D}$  be as in \ref{simplicial}.
  Then
  it follows that
$${\bq}(D_{{\ba}_i})=
  r_i{\bq }({\ba}_i^*)$$  for every $\bq\in |\D|$  and every $i$.
\end{lem}
\begin{pf}
  Let $D_{{\ba}_i}$  be defined by $x_i=0$  on
  $U_{{\ba}_i} \subset V$.
  Then $x_i=y_i^{r_i}$,
  where the equation $y_i=0$  defines  ${\overline D_{{\ba}_i}}$ on
  $T_{\overline N}(\D)$.
  Therefore,
  ${\bq}(D_{{\ba}_i})={\bq}(x_i)=r_i{\bq}(y_i)$
  and the last term equals $r_i {\bq }({\ba}_i^*)$,
  because $y_i$ is the $i$-th coordinate function of
  $T_{\overline N}(\D)={\bC}^r$.
\end{pf}

\begin{prop}
\label{birational canonical}
  Let  $\D' $  be a finite subdivision of  an arbitrary
  finite fan  $\D$  and
  $\varphi :V'=T_N(\D')\to V=T_N(\D)$   be the corresponding
  birational
  morphism.
  Denote the divisor on $V'$ corresponding to  $\bq \in {\D}'(1)$
  by $D'_{\bq}$.
  Assume that  $K_{V}$ is ${\bQ}$-principal.
  Then
  $$K_{V'}=\varphi^* K_V+ \sum _{{\bq}\in {\D}' (1)-{\D} (1)}
  ({\bq } (\varphi^*(\sum_{{\bt}\in {\D} (1)}D_{\bt}))-1)D'_{\bq}.$$

  If $\D$  is as in \ref{simplicial},
  then
  it follows that
   $$K_{V'}=\varphi^*K_V+ \sum _{{\bq}\in {\D}'(1)-{\D}(1)}
   (\sum_{i=1}^r r_i{\bq}({\ba}_i^*)-1)D'_{\bq}.$$

  If moreover $\sigma$  is the positive quadrant in  $N_{\bR}$,
  then
  $$K_{V'}=\varphi^*K_V+ \sum _{{\bq}\in {\D}'(1)-{\D}(1)}
   ({\bq}({\bl})-1)D'_{\bq},$$
  where  ${\bl}=(1,1,\ldots , 1) \in M$.
\end{prop}

\begin{pf}
  For a toric variety,
  the canonical divisor is represented by the sum of all toric
  invariant divisors with coefficient $-1$.
  Therefore

  (1)  $K_{V'} = -\sum _{{\bq}\in \D' (1)-\D (1)}D'_{\bq}-
        \sum_{{\bt}\in \D (1)}D'_{\bt}$.
  \newline
  On the other hand, we represent $K_{V'}$  as

  (2)  $K_{V'} = \varphi^*K_V+
        \sum_{{\bq}\in \D' (1)-\D (1)}m_{\bq}D'_{\bq}$.
  \newline
  Substituting  $K_V=-\sum_{{\bt}\in \D (1)}D_{\bt}$
  into  (2) and comparing (1) and (2),
  we obtain the value of $m_{\bq}$.
  For the second and the third equalities,
   note that  $\bq(\varphi ^* D)=\bq(D)$
  for a $\bQ$-principal divisor $D$ and  apply \ref{dual}.
\end{pf}

\vskip 1truecm

  Now we obtain the characterization of hypersurface singularities
  by means of the
  Newton diagram.
  Parts of the following are stated in
 \cite{HW}, \cite{W3} and \cite{Reid2}.

\begin{cor}
\label{characterization of log-canonical}
  Let  $(X,0) \subset ({\bC}^{n+1}, 0)$  be an isolated
  singularity defined by a polynomial  $f$.
  Denote by  $\Gamma_{+}(f)$  and  $\Gamma(f)$
   Newton's diagram of  $f$  and the union of the compact faces of it,
  respectively.
  Then the following hold:

  (i) if  $(X,0)$  is canonical, then ${\bl}=(1,1,\ldots , 1)\in
  \Gamma_{+}(f)^o$,
  where $\Gamma_{+}(f)^o$is the interior of $\Gamma_{+}(f)$;

  (ii) if $(X,0)$ is  log-canonical, then  ${\bl}\in \Gamma_+(f)$.
\vskip .5truecm
If  $f$  is non-degenerate, the following hold:

  (iii) $(X,0)$  is canonical if and only if  ${\bl}\in
  \Gamma_{+}(f)^o$;

  (iv) $(X,0)$ is non-canonical, log-canonical
  if and only if  ${\bl}\in \Gamma(f)$;

  (v) $(X,0)$ is not log-canonical
  if and only if  ${\bl} \notin
  \Gamma_{+}(f)$.
\end{cor}

\begin{pf}
  Let  $\sigma$  be the positive quadrant in  $N_{\bR}$
  and  $\D$  be the fan consisting of all faces of  $\sigma$.
  For a primitive  ${\bq}\in N\cap {\sigma} $ ,
  take the  subdivision  $\Dq$  of  $\D$
  consisting of all faces  of $\sigma_i=\sum_{j\neq i}\bR _{\geq 0}{\be _j}+
  \bR _{\geq 0}$,  $i=0,\ldots , n$,
  take the normalization
  $\tilde X$  of
  the proper transform
  $\overline{X}\subset T_N(\Dq)$  of  $X$.
  For the composite
  $\psi:\tilde X \to \overline{X}
\stackrel{\varphi|_{\overline{X}}}\longrightarrow X$,
  write the canonical divisor as follows:
  $$K_{\tilde X}=\psi ^* K_X +\sum m_iE_i,$$
  where  $E_i$'s are the exceptional divisors of  $\psi$.
  On the other hand
  $$K_{T_N(\D(\bq))}+\overline{X}=
  \varphi^*(K_{\Cn}+X)+(\bq(\bl)-1-\bq(f))D_{\bq}.$$
  For the statement (i) (resp. (ii)),
  it is sufficient to prove that if  ${\bq}(\bl)-1-{\bq}(f)<0$
  (resp. $<-1$),
  then  $m_i<0$
  (resp.  $m_i<-1$)  for some  $i$.
 Since  $T_{N}(\Dq)$  has at worst ${\bQ}$-factorial
  log-terminal singularities and  $K_X$  is linearly trivial,
  we can apply the following lemma to a Weil divisor  $\overline{X}
  \subset  T_{N}(\Dq)$.
  For the assertion of the case that  $f$  is non-degenerate,
  it is sufficient to prove the opposite implications in (i) and (ii).
  For a non-singular subdivision  $\D'$ of  $\D$,  on whose toric variety
  the proper transform  $X'$  is
  non-singular and intersects transversally each orbit on $T_{N}(\Dq)$
 (for the existence of such  $\D'$, cf. \cite{Kouchinirenco}
  \cite{V}),
  we have:
  $$K_{X'}=\varphi ^*(K_X)+\sum_{\bq \in \D '(1)-\D(1)}
  ({\bq}(\bl)-1-{\bq}(f))D'_{\bq}|_{X'}$$
  by \ref{valuation of function}, \ref{birational canonical}.
  If  $\bl \in \GGf^o$,
  then $\bq(f)<\bq({\bl})$  for all $\bq \in \D '(1)-\D(1)$,
  which implies that  $(X,0)$  is canonical.
  If it is not log-canonical, then there exists  $\bq$ such that
  $\bq(f)>\bq({\bl})$,
  which implies   ${\bl} \notin
  \Gamma_{+}(f)$.
 \end{pf}

\begin{lem}
\label{normalization}
  Let  $Y\subset  Z$  be an irreducible
  Weil divisor on a normal variety  $Z$.
  Suppose Z  admit at worst  $\bQ$-factorial log-terminal singularities.
  Let  $\tau :\tilde Y \to Y$  be the normalization.
  Then:
\newline
  (i)  $Y$  is a Cohen-Macaulay variety;
\newline
  (ii)  $\omega_Y\simeq
  (\omega_Z(Y)\otimes_{{\cal O}_Z}{\cal O}_Y)/{\cal T}$,
  where ${\cal T}$  is the torsion submodule of
  $\omega_Z(Y)\otimes_{{\cal O}_Z}{\cal O}_Y$;
\newline
  (iii)  if  $\omega_Z(Y)\simeq {\cal O}_Z(-aD)$  $a\geq 1$  (resp.
  $a>1$ )  for an effective divisor  $D\subset Z$  such that
  $\phi \neq  D\cap Y \neq Y$,
   then we have the canonical isomorphism
  $\omega _{\tilde Y} \simeq {\cal O}_{\tilde Y}(-\sum_{i=1}^kb_iE_i)$
  with  $b_i\geq 1$  for divisors  $E_i$ $(i=1,\ldots , k)$ such that
  $\cup_{i=1}^kE_i\supset \tau^{-1}(D)$
  (resp. in addition  $b_i>1$  for some $i$  such that
  $E_i\subset Supp(\tau^{-1}D)$).
\end{lem}
\begin{pf}
  First, one can prove that every effective Weil divisor on $Z$  is a
  Cohen-Macaulay variety in the same way as in 0.5 of  \cite{IF},
  because the covering constructed as in \cite{IF} has at worst
  rational singularities in the present case too.
  For the proof of  (ii),
  take the exact sequence:
$${\cal Hom}_{{\cal O}_Z}({\cal O}_Z ,\omega_Z)\to
  {\cal Hom}_{{\cal O}_Z}({\cal O}_Z(-Y) ,\omega_Z)\to
  {\cal Ext}^1_{{\cal O}_Z}({\cal O}_Y,\omega_Z)\to
  {\cal Ext}^1_{{\cal O}_Z}({\cal O}_Z,\omega_Z)=0.
$$
  Here ${\cal Ext}^1_{{\cal O}_Z}({\cal O}_Y,\omega_Z)=\omega_Y$,
  because  $Z$  is a Cohen-Macaulay variety.
  So  $\omega_Y$  is the image of
   ${\cal Hom}_{{\cal O}_Z}({\cal O}_Z(-Y) ,\omega_Z)=\omega_Z(Y)$,
  and therefore also the image of $\omega_Z(Y)\otimes {\cal O}_Y$
  which is isomorphic to  $\omega_Y$ on general points of  $Y$.
Since  $\omega_Y$  is torsion-free,
  it must be isomorphic to
  $(\omega_Z(Y)\otimes_{{\cal O}_Z}{\cal O}_Y)/{\cal T}$
  as desired in (ii).
  Now,
  since  $\tau$  is finite,
  we have the inclusion $\tau_*\omega_{\tilde Y}\hookrightarrow
  \omega_Y$.
  By (ii) and the assumption of (iii),
  $\omega_Y$  is isomorphic to the defining ideal  ${\cal J}$  of
  a subscheme $aD\cap Y $ of $ Y$.
  Therefore  $\omega_{\tilde Y}\simeq {\cal O}_{\tilde Y}(-\sum b_iE_i)$,
    $b_i> 0$  for divisors  $E_i$ $(i=1,\ldots , k)$ such that
  $\cup_{i=1}^kE_i\supset \tau^{-1}(D)$.
  Next, assume  $a>1$.
  For the assertion,
  we may replace  $Y\subset Z$  with a small neighbourhood of a general
  point on  $D\cap Y$.
  So we may assume that all $E_i$ are over $Supp(D|_Y)$ and $D|_Y$  is
  irreducible.
  If there is no $E_i\subset Supp(D)$  such that  $b_i>1$,
  then  $\tau_*\omega_{\tilde Y}=\tau_*{\cal O}_{\tilde Y}(-\sum E_i)$
  is a reduced  ${\cal O}_Y$-ideal whose locus has the support on
  $D\cap Y$.
  On the other hand,  ${\cal J}$  also has the locus with the support on
  $D\cap Y$,
  therefore they coincide.
  By this equality  $\tau_*\omega_{\tilde Y}=\omega_Y$, it follows that
$$\tau_*{\cal O}_{\tilde Y}={\cal Hom}_{\tau_*{\cal O}_{\tilde Y}}
(\tau_*{\omega}_{\tilde Y},\tau_*{\omega}_{\tilde Y})
\subset
{\cal Hom}_{{\cal O}_{ Y}}
({\omega}_{ Y},{\omega}_{Y})={\cal O}_Y,$$
  where the left and right equalities follow from the fact that
  $\tilde Y$
  and $Y$  satisfy  $S_2$-condition.
  Now it follows that  $\tau: \tilde Y \simeq Y$  is normal,
  which induces
the contradiction to  $a>1$.
\end{pf}

\begin{say}
\label{kappa}
  For a normal isolated singularity  $(X,x)$, we define an invariant
$\kappa_{\delta}(X,x)$(\cite{Iasym} )
by the growth order of the plurigenera  $\delta_m$ $(m \in {\bN})$
(\cite{W}).

In general, n-dimensional, normal, isolated singularitieis  $(X,x)$  are
classified by the invariant $\kappa_{\delta}$ into (n+1)-classes:
$\kappa(X,x) = -\infty,\ 0, \ 1, \ \ldots , \ n-2, \ n$ (skipping $n-1$
curiously)
(\cite{Iasym}).
For hypersurface singularities,
the classes are only three :
$\kappa(X,x) = -\infty,\ 0,\ n$ (\cite{TW}).
A hypersurface singularity with  $\kappa_{\delta}(X,x)=-\infty$
  (resp.$ =0,\ n$)
  is equivalent to the fact that $(X,x)$  is canonical (resp.
  non-canonical-log-canonical, not log-canonical) (cf. \cite{IGor}).
  Therefore  \ref{characterization of log-canonical} also gives
  the  combinatoric characterization of non-degenerate hypersurface
  singularities' classes by $\kappa_{\delta}$
\end{say}


\section{The weights which give the canonical modification}

\begin{say}
\label{quadrant}
  Under the notation in \ref{simplicial},
  put $r=n+1$  for $n\geq 2$  and number the elements of the basis
  $\{\be_i\}$  from $i=0$  to $i=n$.
  Let   $\sigma=\sum _{i=0}^n{\bR}_{\geq 0} \be_i$  be
the positive quadrant in  $N_{\bR}$,
  and  $\D=<\sigma>$  be the fan consisting of all faces of  $\sigma$.
  Denote  Newton's diagram of a polynomial  $f\in {\bC}[x_0,\ldots ,x_n]$
  and the union of its compact faces
  by $\GGf$ and by  $\Gf$ respectively.
\end{say}
\begin{defn}
\label{essential cone}
  For a polynomial  $f \in {\bC}[x_0,\ldots ,x_n]$,
  we define the essential cone as follows:
   $${\Cf}:=\{ \bq \in \sigma \subset N_{\bR}|\bq(f)-\bq(\bl)\geq 0\}.$$

\end{defn}

\begin{rem}
\label{rem of essential cone}
  (i)  It is clear that if $\bl \in \GGf ^o$,
  then $\Cf =\{ 0 \}$.

  (ii)  If $\bl \notin \GGf^o$, the  essential cone  $C_{\bf 1}(f)$  is
actually
  the cone spanned by  $\gamma_1^{\perp } ,\ldots , \gamma_r^{\perp } $,
  where each  $\gamma_i$  is an n-dimensional face of
  $\Gamma_+(x_0\cdots x_n+f)$
  which contains ${\bl}$.
  Let  $X$  be the divisor in  $\bC ^{n+1}=T_N(\D)$  defined by  $f=0$.
  If  $X$  has an isolated singularity at the origin  $0 \in \Cn$,
  then every vector  $\bq \in \Cf-\{ 0\}$  has positive coordinates
  $q_j$ for  $j=0,1,\ldots , n$,
  otherwise at least one  $\gamma_i$ is parallel to one of the coordinate
  axes which causes a contradiction to the isolatedness of the singularity
$(X,0)$.

   (iii)  in  Def 3.3 of \cite{IGor} ,
   we have the notion of an essential divisor of a resolution of a Gorenstein
   singularity.  Every 1-dimensional cone in the essential cone  in
\ref{essential
   cone} gives a component of the essential divisor in some resolution.
\end{rem}
\begin{defn}
\label{order}
  (1)  Let  $C$  be a cone in $\sigma \subset N_{\bR}$.
  For  $\bp =(p_0,\ldots , p_n),\ \ \bq = (q_0,\ldots ,q_n) \in C$,
  we define
  $\bp\leq \bq$  if $p_i\leq q_i$  for every $i=0,\ldots ,n$.

  We say that  a primitive element $\bp \in C\cap N-\{ 0\}$  is
  absolutely minimal,
  if $\bp\leq \bq$  for every primitive element
  $\bq \in   C\cap N-\{ 0\}$.

  (2) For  $\bp,\  \bq \in \Cf$,
  we define
  $\bp\leq _f\bq$ ,
  if  $p_i/(\bp(f)-\bp(\bl)+1 )\leq q_i/(\bq(f)-\bq(\bl)+1)$
  for every  $i=0,\ldots ,n$.
  We define another order  $\prec _f$ as follows:
  $\bp \prec _f \bq$ if
  $p_i/(\bp(f) )\leq q_i/(\bq(f))$
  for every $i=0,\ldots ,n$.
  We say that  a primitive element $\bp \in \Cf\cap N-\{ 0\}$
  is  $f$-minimal,
  if for every primitive element
  $\bq \in   \Cf\cap N-\{ 0\}$,
  either $\bp\leq _f\bq$
  or $\bp\prec _f \bq$ and $\bq$ belongs to the interior of an
$n+1$-dimensional
   cone of  $\D(\bp)$,
   where the fan $\D(\bp)$ consists of all faces  of
  $\sigma_i=\sum_{j\neq i}\bR_{\geq 0}\be_j +\bR_{\geq 0}\bp
  \subset N_{\bR}$,  $i=0,\ldots , n$.
\end{defn}
\begin{say}
\label{notation star}
  For a primitive vector  $\bp \in \sigma\cap N-\{ 0\}$,
  we have the star-shaped decomposition  $\D(\bp)$
  by adding the ray  $\bR_{\geq 0}{\bp}$ as in the  definition above.
  We denote the fan of all faces of $\sigma_i$ by $\D_i$.
  Denote the proper transform of  $X=\{f=0\}$  on
  $T_N(\D(\bp))$ by $X(\bp)$.
 The induced morphisms
  $\varphi : T_N(\D(\bp))\to T_N(\D)$,
  $\varphi ':X(\bp) \to X$  are called weighted blow-ups with weight
 $\bp$,
  or simply $\bp$-blow-ups of $\Cn$ and $X$  respectively.
  Let  $U_i$  be the invariant open subset
  $U_{\sigma _i} \simeq T_N(\D_i)$  of  $T_N(\D(\bp))$,
  and  $\varphi_i:U_i \to \Cn$  be the restriction of  $\varphi$ onto
  $U_i$.
  Denote  $X(\bp)\cap U_i$  by  $X_i$.
\end{say}

\begin{prop}
\label{prop star}
  Under the notation in \ref{notation star},
  let  $\psi:\tilde{U_i}\to U_i$  be the birational morphism corresponding
  to a finite subdivision $\Sigma_i$  of  $\D_i$.
  Denote the proper transform of  $X_i$  by  $\tilde{X_i}$,
  then  $$K_{\tilde{U_i}}+\tilde{X_i}=
  \psi_i^*(K_{U_i}+X_i)
  +\sum_{\bq\in \Sigma_i(1)-\D_i(1)}(\dfrac{q_i}{p_i}
  (\bp (f)-\bp(\bl)+1)-(\bq (f)-\bq(\bl)+1))D_{\bq}.$$
\end{prop}
\begin{pf}

 We can assume that  $i=0$  without the loss of generality.
  Let  $\{\ba_j\}_{j=0}^n$  be  $\{\bp, \be _1, \ldots ,\be _n\}$.
  By \ref{valuation of function} and \ref{birational canonical},
  it is sufficient to prove:
  $$(1)\ \ \ \ \ \
   \sum_{j=0}^nr_j\bq(\ba^*_i)-\bq(\psi^*X_0)-1=
  \dfrac{q_0}{p_0}(\bp (f)-\bp(\bl)
  +1)-(\bq (f)-\bq(\bl)+1))D_{\bq}.$$
  First, we can see that  $r_j=1$  for every  $j$.
  In fact,
  the quotient map  $\pi:\Cn=T_{\overline{N}}(\D_0)\to
  T_{N}(\D_0)=U_0$   is defined by the action of the cyclic group
  generated by

$$\pmatrix
  \epsilon  & 0               & \ldots         & \ldots & 0\\
  0         & \epsilon^{-p_1} & 0              &\ldots  & 0\\
  0         &         0       &\epsilon^{-p_2} & 0      &\vdots\\
   \vdots   & \vdots          &\vdots          &\ddots  &\vdots\\
  0         & 0               &    0           &\hdots  &\epsilon^{-p_n}
\endpmatrix ,$$

  where  $\epsilon$  is a primitive $p_0$-th root of unity.
  Here it is easy to check that  $\pi$  is etale in codimension one.
  Next, since  $\ba_0^*=\bp^*=(1/p_0,0,\ldots ,0) $  and
  $\ba_j^*=\be_j^*=(-p_j/p_0,0,\ldots ,0,1,0,\ldots ,0)$
  (j-th entry is 1) for $1\leq j \leq n$,
  one obtains $\bq (\sum \ba_j^*)=
  (1-p_1-\ldots -p_n)q_0/p_0+(q_1+\ldots +q_n)$.
  On the other hand,
  since  $\varphi_0^*(X)=X_0+\bp (f)D_{\bp}$
  by \ref{valuation of function},
  it follows that
$$(2)\ \ \ \ \ \
\bq(\psi^*X_0)=\bq(\psi^* \varphi_0^*(X))-
  \bq(\psi^*(\bp(f)D_{\bp}))
=\bq(f)- \bp(f)\bq(\bp ^*)
=\bq(f)- \bp(f)q_0/p_0.$$
  By substituting them into the left hand side of (1)  we obtain the
  equality (1).
\end{pf}
\vskip 1truecm
\begin{lem}
\label{normality}
 Let  $Y\subset Z$  be an irreducible  Weil divisor on a variety
   $Z$. Assume that  $Z$ admits at worst $\bQ$-factorial log-terminal
singularities.
  Let  $\psi:\tilde Y \to Y$  be a resolution of singularities on  $Y$.
  Assume $K_{\tilde Y}=\psi^*((K_Z+Y)|_{Y})+\sum_im_i E_i$
  with  $m_i>-1$  for all $i$,
  where $E_i$'s are the exceptional divisors of  $\psi$.

   Then  $Y$  is normal,
  and $Y$  has at worst log-terminal singularities.

  In particular,
if $m_i\geq 0$ for all $i$,
then  $Y$  has at worst canonical singularities.
\end{lem}

\begin{pf}
First  $Y$  is a Cohen-Macaulay variety as in \ref{normalization}.
  Therefore it is sufficient to prove that $codim _YSing(Y)\geq 2$
  by Serre's criterion.
  Assume that  $y\in Y $  is a general point of a component of $Sing( Y)$
  of codimension one.
 By replacing $Y$  with a small neighbourhood of $y$,
  we may assume that  $\psi$  is the normalization.

Claim that $\psi_*\omega_{\tilde Y}=\omega _Y$.
  The inclusion $\subset$  is trivial.
  For the proof of the opposite inclusion,
  take an arbitrary  $\theta\in \omega_Y$.
  Then  $\theta^r\in \omega_Z^{[r]}(rY)\otimes {\cal O}_Y$
for such  $r$ that $\omega_Z^{[r]}(rY)$  is invertible,
  because  $\omega_Y=\omega_Z(Y)\otimes{\cal O}_Y/{\cal T}$
  by (ii) of \ref{normalization}.
  By the assumption of the lemma, one obtains:
$$\theta^r\in \omega_Z^{[r]}(rY)\otimes {\cal O}_Y\subset \psi_*\omega
_{\tilde Y}^r(-\sum rm_iE_i).$$
  Hence for the valuation  $\nu_i$  at each $E_i$,
  $r\nu_i(\psi^*\theta)=\nu_i(\psi^*\theta^r)\geq rm_i >-r$.
  Therefore $\nu_i(\psi^*\theta)\geq 0$  for every  $E_i$,
which means  that $\psi^*\theta\in \omega_{\tilde Y}$  as claimed.
  By the same argument as in the proof of (iii) in \ref{normalization},
  it follows that  $Y$  is normal.
  One can see also that $Y$  is  $\bQ$-Gorenstein,
  because $\omega_Y^{[r]}=\omega_Z^{[r]}(rY)\otimes {\cal O}_Y$
  is invertible for $r$ above.
\end{pf}
\vskip 1truecm
\begin{thm}
\label{main theorem}
  Let  $(X,0)\subset(\bC ^{n+1},0)$  be an isolated singularity defined
  by a polynomial  $f\in \bC[x_0,\ldots , x_n]$.
  For a primitive integral vector  $\bp =(p_0,\ldots ,p_n)$
  such that  $p_i>0$ for all $i$,
\newline
  if

  (i)  $\bp$  is the canonical weight, i.e.,
   $\bp$-blow-up  $\varphi :X(\bp)\to X$  is
      the canonical modification,
\newline
  then

  (ii) $\bp$  is  $f$-minimal in $\Cf \cap N-\{0\}$.
\vskip .5truecm

  Suppose  $f$  is non-degenerate,
  then the converse  (ii)$\Rightarrow$(i)   also holds.

\end{thm}
\begin{pf}
  We use the notaion in \ref{notation star}.
  If the $\bp$-blow-up $\varphi':{X(\bp)}\to X$  is
  the canonical modification,
  then it follows that  $\bp \in \Cf$;
  otherwise,
  $K_{X(\bp)}=\varphi^*K_X+(\bp(\bl)-1-\bp(f))D_{\bp}|_{X(\bp)}$  with
  $\bp(\bl)-1-\bp(f)\geq 0$,
  which shows a contradiction
  that  $(X,0)$ itself is a canonical singularity .
  Let us prove that  $\bp$ is  $f$-minimal.
  If there exists  $\bq \in \Cf \cap N-\{0\}$
  such that  $\bq \not\geq _f \bp$,
  then $\min\{\dfrac{q_i/(\bq(f)-\bq(\bl)+1)}{p_i/(\bp(f)-\bp(\bl)+1 )}
  | i=0,\ldots , n\}
      <1$.
  Let  $i=0$  attain the minimal value,
  then it follows that $\bq\in \sigma_0$,
  because $\bq$  is represented as
  $\dfrac{q_0}{p_0}\bp+ \sum _{i=1}^n(q_i-\dfrac{q_0}{p_0}p_i)\be_i$
  and its coefficients are all non-negative.
  Taking the star-shaped subdivision  $\D_0(\bq)$  of  $\D_0$
  by adding a ray  $\bR_{\geq 0}\bq$,
  we have a birational morphism
  $\psi: \tilde{U_0}:=T_N(\D_0(\bq)) \to U_0=T_N(\D_0)$.
  Denote the proper transform of $X_0$  by  $\tilde{X_0}$,
  then, by \ref{prop star}
  $$(K_{\tilde{U_0}}+{\tilde{X_0}})|_{\tilde{X_0}}=
  \psi^*(K_{X_0})
  +(\dfrac{q_0}{p_0}
  (\bp (f)-\bp(\bl)+1)-(\bq (f)-\bq(\bl)+1))D_{\bq}|_{\tilde{X_0}}.$$
  It follows that ${\tilde{X_0}}$  is normal.
In fact, for a resolution $\lambda:{\tilde{\tilde{X_0}}} \to {\tilde{X_0}}$,
  denote
$K_{\tilde{\tilde{X_0}}}=\lambda^*((K_{\tilde{U_0}}+{\tilde{X_0}})|_{\tilde{
X_0}})
  +\sum m_iE_i$,
  then $K_{\tilde{\tilde{X_0}}}=\lambda^*\psi^*K_{X_0}+\sum (n_i+m_i)E_i$
  where $n_i$ is the coefficient of $E_i$  in  $\lambda^*((\dfrac{q_0}{p_0}
  (\bp (f)-\bp(\bl)+1)-(\bq (f)-\bq(\bl)+1))D_{\bq}|_{\tilde{X_0}})$
  which is  non-positive by the negativity of
  $\dfrac{q_0}{p_0}
  (\bp (f)-\bp(\bl)+1)-(\bq (f)-\bq(\bl)+1)$;
  therefore  if  $m_i<0$ for some  $i$,
  then it contradicts to the fact that  $X_0$  is  canonical;
  since $m_i \geq 0$ for all $i$,
  by Lemma \ref{normality}  $\tilde{X_0}$  is normal.
  Now we obtain the partial resolution $\psi: \tilde{X_0} \to X_0$
  with $K_{\tilde{X_0}}=
   \psi^*(K_{X_0})
  +(\dfrac{q_0}{p_0}
  (\bp (f)-\bp(\bl)+1)-(\bq (f)-\bq(\bl)+1))D_{\bq}|_{\tilde{X_0}}$.

  If $D_{\bq}\cap \tilde{X_0}\neq \phi$,
  the coefficient of $D_{\bq}|_{\tilde{X_0}}$  is negative
  by the definition of  $\bq$,
  which contradicts the hypothesis that $X_0$  has at
  worst canonical singularities.
  Therefore,  $D_{\bq}\cap \tilde X_0=\phi$
  which happens if and only if  $\psi(D_{\bq})$  is a point away from
  $X(\bp)$, because $X(\bp)$  is ample on  $D_{\bp}$.
  It implies
  $\bq(\psi ^*X_0)=\bq(f)-\dfrac{q_0}{p_0}\bp(f)=0$
  (c.f. (2) in the proof of \ref{prop star})
   and  $\bq$ belongs to the interior of an $n+1$-dimensional
   cone of  $\D(\bp)$.
  By the minimality of  $q_0/p_0$,
  we obtain  $\bp \prec _f \bq$.

  Next suppose that  $f$  is non-degenerate and  $\bp$  is
  $f$-minimal in $\Cf\cap N-\{0\}$.
  Take a non-singular subdivisioin  $\D' $ of $\D(\bp)$
  such that the restriction of the corresponding morphism
  $\psi:T_N(\D')\to T_N(\D(\bp))$  onto the proper transform
  $X(\D')$  of  $X(\bp)$  gives a resolution of $X(\bp)$
  such that  every intersection of $X(\D')$ and an orbit is transversal.
  Let  $\D'_i$  be the subdivision of  $\D_i$  which is in $\D'$,
  then  $T_N(\D')$  is covered by  $T_N(\D'_i)$'s
  and the restriction $\psi_i: T_N(\D'_i)\to U_i=T_N(\D_i)$
   of $\psi$  gives a resolution  $X(\D'_i):=X(\D')\cap T_N(\D'_i)\to X_i$
  of each  $X_i$.

 Represent the canonical divisor on $X(\D'_i)$ by
$$
  K_{X(\D'_i)}= \psi_i^*(K_{U_i}+X_i)|_{X(\D'_i)}+\sum_{\bq\in\D'_i(1)-
 \D_i(1)}
  m_{\bq}(D_{\bq}\cap {X(\D'_i)})_{red}$$
  If  $D_{\bq}\cap {X(\D'_i)}\neq \phi$,
  the both intersect generically transversally each other by the
construction of
  $\psi$.
  Therefore
  $m_{\bq}
  =\dfrac{q_i}{p_i}
  (\bp (f)-\bp(\bl)+1)-(\bq (f)-\bq(\bl)+1)$.
  If  $\bq\notin \Cf$,
  then  $\bq (f)-\bq(\bl)+1\leq 0$
  and therefore  $m_{\bq}>0$.
  If $\bq\in \Cf$ and $D_{\bq}\cap {X(\D'_i)}\neq \phi$,
  then ${\bp}\not\prec _f{\bq}$.
  In fact,
  if   ${\bp} \prec _f{\bq}$,
  then ${\bq}  $  is in the interior of  ${\D_i}$ for some  $i$,
  let it be  ${\D_0}$,which implies  $\psi_0(D_{\bq})$  is a point.
  We also have that  ${\bq}(f)-\dfrac{q_0}{p_0}{\bp}(f)\leq 0$.
  On the other hand,
  ${\bq}(\psi_0^*X_0)={\bq}(f)-\dfrac{q_0}{p_0}{\bp}(f)\geq 0$.
  Therefore ${\bq}(\psi_0^*X_0)=0$, which shows that  $X(\D_0')\cap
D_{\bq}=\phi$.
  Thus it follows that $m_{\bq}\geq 0$ by the absolute $f$-minimality of
$\bp$.
  Now, by Lemma \ref{normality},
  it follows that  $X_i$'s have at worst canonical singularities.

  The $\varphi$-ampleness of $K_{X(\bp)}$  follows from
  the $\varphi$-ampleness of $-D_{\bp}|_{X(\bp)}$ and
$K_{X(\bp)}=\varphi^*K_X+(\bp(\bl)-\bp(f)-1)D_{\bp}|_{X(\bp)}$,
  where the coefficient of  $D_{\bp}$ is negative.
\end{pf}

\begin{cor}
\label{log-canonical}
  Let  $(X,0)\subset(\Cn,0)$  be an isolated, non-canonical,
  log-canonical singularity defined by a polynomial  $f\in \bC[x_0,
\ldots , x_n]$.
  For a primitive integral vector  $\bp=(p_0,\ldots ,p_n)$
  such that  $p_i> 0$  for all $i$,
\newline
if

  (i)  $\bp$-blow-up  $\varphi :X(\bp)\to X$  is
      the canonical modification,
\newline
then

  (ii) $\bp$  is absolutely minimal in $\Cf \cap N-\{0\}$.
\vskip .5truecm

  Suppose  $f$  is non-degenerate,
  then the converse  (ii)$\Rightarrow$(i)  holds too.

\end{cor}

\begin{pf}
Since  $(X,0)$  is log-canonical,
  $\bl\in\GGf$ by \ref{characterization of log-canonical}.
  Then  $\bq(f)=\bq(\bl)$, for $\bq\in \Cf$.
  First,
two distinct primitive $\bp$, $\bq \in \Cf$,
  neither $\bp \prec _f  \bq$ nor $\bq \prec _f  \bp$ hold.
  In fact,
  if  $\bp \prec _f  \bq$, then  $p_i/\bp(\bl)\leq q_i/\bq(\bl)$
  for every  $i$,  and moreover the equality holds for every $i$,
  because,
  by summing all these inequalities,
  we obtain $1\leq 1$.
  Hence $\bp$  must coincide with $\bq$.
  On the other hand,
  it is clear that  $\geq_f$ is equivalent to $\geq$ and therefore

  $f$-minimal is equivalent to absolutely minimal.
\end{pf}.

\begin{exmp}
  (Tomari)
  It is possible for a singularity to have more than one canonical weights.
  In fact,
  let  $X\subset \bC ^3$  be defined by $x_0^k+x_1^{k+1} +x_2^{k+1}=0$
$(k\geq 3)$,
  then the weights  $(1,1,1)$ and $(k+1, k, k)$ are both canonical weights.
\end{exmp}

\begin{say}
\label{various cor}
  In the rest of this section a singularity  $(X,0)\subset(\Cn,0)$
  is assumed to be  a non-canonical, log-canonical singularity
  defined by a  polynomial  $f$.

  In some cases,
  one can easily see the existence of the absolutely minimal vector,
  therefore one also sees the existence of the canonical modification  for
these cases.

  For 2-dimensional case,
  singularities as above are either $\tilde{E_6}$ or $\tilde{E_7}$ or
  $\tilde{E_8}$ or defined by equations of type:
  $x_0x_1x_2+x_0^p+x_1^q+x_2^r=0$ with  $\dfrac{1}{p}+
  \dfrac{1}{q}+\dfrac{1}{r}<1$
   by suitable coordinates transformations.
  It is well known that there exist the canonical weights
  for  singularities defined by these equations.
  The following corollary shows that one need not take a coordinate
  transformation for 2-dimensional non-canonical,
  log-canonical singularities to
  admit the canonical weight.
\end{say}
\begin{cor}
\label{surface}
  If  $n=2$ and  $f$  is non-degenerate,
  there exists the absolutely minimal vector in  $\Cf$  for every $f$
  as in  \ref{various cor}.
  And the vector is either  $(1,1,1)$  or $(3,2,1)$  or $(2,1,1)$.
 \end{cor}
\begin{pf}
  By the direct calculation, one can find the absolutely minimal vector
  in $\Cf$ for each  $f$.
\end{pf}
  The next one was  proved by Tomari
  under a more general situation.
\begin{cor}
\label{Tomari's cor}
 (\cite{T})
  If ${\bl}$ is in the interior of an n-dimensional face $\gamma$
  of  $\Gamma(f)$,
  which is equivalent to that the singularity  $(X,0)$  is of Hodge type
  $(0,n-1)$ (for the definition, cf. \cite{IGor}),
  then the primitive vector ${\bp}$ generating  $\gamma^{\perp } $
  gives the canonical modification  $X({\bp})\to X$.

\end{cor}
\begin{pf}

  Since  $ C_{\bf 1}(f)$ is of one dimension,
  the primitive vector on it is clearly absolutely minimal.
  It completes the proof of the non-degenerate case.
  If $f$  is degenerate,
  there may not exist toric embedded resolution.
  But  taking a resolution  $\psi:Y\to X(\bp)$,
  $\varphi\psi$ is a  resolution of
  a  log-canonical singularity of type  $(0, n-1)$,
  which yields that $K_Y=\psi^*\varphi^*(K_X)+\sum_i m_iE_i$
  with the only one  negative $m_i$.
  By substituting
  $(K_{T_N(\D(\bp))}+X(\bp))|_{X(\bp)}=
  \varphi^*(K_X)-D_{\bp}|_{X(\bp)}$ into the equality above,
  we can see that
   the pair $X(\bp)\subset T_N(\D(\bp))$ satisfy the conditions of
  \ref{normality}.
\end{pf}
\begin{cor}
  If  a non-degenerate polynomial $f$  is represented as  $x_0\cdots x_n
+h(x_0,\ldots ,x_n)$,
  where $deg\ h\geq n+1$,
  then the blow-up  by the maximal ideal of the origin is
  the canonical modification.
\end{cor}
\begin{pf}
  Since  $\GGf$  is in the domain
  $\{\ba\in M_{\bR}|a_0+\cdots +a_n \geq n+1\}$,
  it follows that  $(1,1,\ldots ,1)\in \Cf$ and clearly this is absolutely
  minimal.
\end{pf}

\begin{cor}
\label{quasi-reduced}
  If every vector $\ba\in \Gf\cap M$  is quasi-reduced
  (i.e. $\ba=(a_0,\ldots , a_n)$  satisfies that $0\leq a_i\leq 1$
  except for at most one $i$),
  then there exists the absolutely minimal vector $\bp$  in $\Cf$.
\end{cor}
\begin{pf}
  A positive vector $\bq\in N$ belongs to  $\Cf$,
  if and only if  $\bq(\ba)\geq\bq(\bl)$  for all $\ba=(a_0,\ldots ,a_n)
  \in \Gf\cap M$.
  These inequalities are equivalent to the inequalities of the following
type:
  $(a_i-1)q_i\geq\sum_{j\in \Lambda(\ba)}q_j$,
  where  $a_i\geq 2$  and $\Lambda(\ba)$  is the suitable subset of
  $\{0,\ldots ,n\}$ such that  $i\notin \Lambda(\ba)$.
  Let  $\bp=(p_0, \ldots ,p_n)$ and $\bq=(q_0,\ldots ,q_n)$
  belong to $\Cf$.
  Define  $\br=(r_0,\ldots ,r_n)$  by  $r_i=min\{p_i,q_i\}$.
  We show that  $\br\in \Cf$.
  For  $\ba\in\Gf\cap M$,
  let  $a_i\geq2$.
  We can assume that $r_i=p_i$ by the definition of $\br$.
  Then $(a_i-1)r_i=(a_i-1)p_i\geq
  \sum_{j\in\Lambda(\ba)}p_j\geq\sum_{j\in\Lambda(\ba)}r_j$.
  Hence $\br$ also satisfies  $\br(\ba)\geq\br(\bl)$
  for all $\ba\in \Gf\cap M$.
\end{pf}

\begin{exmp}
\label{typeT}
  We say that   $X$ is of type $T_{\ba}$,
  if  it is defined by $f=x_0\cdots x_n+\sum x_i^{a_i}$  for
  $\ba=(a_0,\ldots, a_n)$,
  where $\sum_i1/a_i<1$.
  Then  $f$  satisfies the condition of \ref{quasi-reduced}
  and therefore $X$  has a weight which gives the canonical modification.
  The summary paper \cite{IB} contains the table of 3-dimensional
  $T_{p,q,r,s}$-singularities  $(X,0)$  with the absolutely minimal vectors
  $\bp$.
  All those weights are in the weights of 95-simple K3-singularities
  listed in \cite{Y} which is bijective to the list of \cite{F}.
  And therefore  $T_{p,q,r,s}$-singularities have the same plurigenera
  $\{\gamma_m\}$ with those of corresponding simple K3-singularities
(cf. \ref{deform 1}).
\end{exmp}
\vskip 1truecm

\section{Deformations and the simultaneous canonical modifications}

\begin{defn}
\label{simul cano}
  Let  $\pi :({\cal X},x) \to (C, 0)$  be a flat morphism
  over a non-singular curve  $C$.
  A partial resolution  $\Phi :{\cal Y}\to {\cal X}$  is called
  the simultaneous canonical modification,
  if the restriction
  $\Phi_t :{\cal Y}_t\to {\cal X}_t$ is the canonical modification for
  every  $t \in C$,
  where ${\cal X}_t=\pi ^{-1}(t)$ and  ${\cal Y}_t=\Phi ^{-1}({\cal X}_t)$.
\end{defn}

\begin{prop}
\label{deform 1}
  Let  $(X,0)\subset (\Cn, 0)$  be an isolated, non-canonical,
log-canonical singularity defined by a polynomial  $f$.
  Assume that  $X(\bp)\to X$  is the canonical modification for a
  positive integral vector $\bp$.
  Let  $\{F_t\}_{t\in C}$  be a deformation of  $f=F_0$ over a non-singular
  curve $C$
  such that $F_t$'s ($t\neq 0$) are non-degenerate and  Newton's
  diagrams $\Gamma_+(F_t)$  sit in the  halfspace $\bl+\bp^{\vee}$
  of  $M_{\bR}$.
  Then
  the flat family $\pi:({\cal X},0)\to (C,0)$  defined by
  $\{F_t\}_{t\in C}$
  admits the simultaneous canonical modification and
  $\gamma_m({\cal X}_t,0)$  is constant in $t\in C$  for every
  $m\in \bN$.
\end{prop}

\begin{pf}
  By the assumption of $\{F_t\}_{t\in C}$,
  $\bp(\bl)\leq \bp(F_t)\leq \bp(f)$.
  Therefore $\bp \in C_{\bl}(F_t) \subset \Cf$.
  Since  $\bp$  is absolutely minimal in $\Cf$,
  it is absolutely minimal in $C_{\bl}(F_t)$ for every  $t\in C$,
  which yields that  $\bp$  is the canonical weight for the singularities
  defined by  $F_t=0$.
    On the other hand,
  since  $(X,0)$  is log-canonical,  $\bl\in \GGf$  by
  \ref{characterization of log-canonical}.
  Hence $\bp(f)=\bp(\bl)$  which implies also  $\bp(\bl)=\bp(F_t)$.
  Take the morphism  $\Phi:=\varphi\times id_C: {T_N(\D(\bp))}
  \times C \to {\bC}^{n+1} \times C$,
  where  $\varphi : T_N(\D(\bp)) \to {\bC}^{n+1}$
  is the $\bp$-blow-up.
  Denote the proper transform of  ${\cal X}$  in  $T_N(\D(\bp))
  \times C $  by  ${\cal Y}$,
  then  $\Phi^*{\cal X}={\cal Y}+\bp(F_t)(D_{\bp}\times C)$ for a
  general $t\in C$,
  where  $D_{\bp}$  is the corresponding divisor to  $\bp$ on
  $T_N(\D(\bp))$.
 Since  $\bp(F_t)$  is constant for all $t\in C$,
  ${\cal Y}_t$'s are all irreducible, and therefore
  these turn out to be the $\bp$-blow-ups of
  ${\cal X}_t$,
which shows that ${\cal Y}\to {\cal X}$  is the simultaneous canonical
  modification.
  By Proposition 7 of  \cite{Stevens},
${\cal Y}$  admits at worst canonical singularities, and,
  on the other hand,
$K_{{\cal Y}/{\cal X}}=-(D_{\bp}\times C)|_{\cal Y}$  is
  $\Phi$-ample,
  which means that  ${\cal Y}\to {\cal X}$  is
  the canonical modification
  of  ${\cal X}$.
  Thus  $\pi$ turns out to be an (FG)-deformation in terms of \cite{Isiml}.
  By 1.11 of \cite{Isiml},
  it follows that  $\gamma_m({\cal X}_t, 0)$  is constant for all $t\in C$.
\end{pf}
\begin{prop}
\label{deform 2}
  Let  $(X,0)\subset (\Cn, 0)$  be an isolated, non-canonical,
log-canonical singularity defined by a polynomial  $f$.
  Assume that  $X(\bp)\to X$  is the canonical modification for a
 positive integral vector $\bp$.
  If  $\bp/\sum_ip_i$ is the weight of a weighted-homogeneous
  polynomial  defining an isolated singularity
  at the origin,
  then  $\gamma_m(X,0)=\gamma_m(Y,0)$,
  where $(Y,0)\subset (\Cn, 0)$ is defined by a non-degenerate
  weighted-homogeneous polynomial $g$ with the weight $\bp/\sum_ip_i$.
  Moreover there exists a flat deformation  $\pi:({\cal X}, 0)\to
  (C,0)$  of $(X,0)=({\cal X}_0, 0)$  over a non-singular curve $C$
  with $({\cal X}_{\tau}, 0)\simeq (Y,0)$  for some  $\tau\in C$
  such that  $\pi$  admits the simultaneous canonical modification.
\end{prop}

\begin{pf}
   Let  $F_t$ be  $(1-t)f+tg$  for  $t\in \bC$.
  Then, taking a suitable open subset  $C\subset \bC$ with $0, 1 \in C$,
  it follows that $F_t$ $(t\neq 0 )$ defines a non-canonical, log-canonical
  singularity of type  $(0, n-1)$, because it is a small deformation
  of such a singularity  $\{g=0\}$ and $\bl\in \Gamma_+(F_t)$  (4.4 of
\cite{def},
  2.2 of \cite{Isiml} and \ref{characterization of log-canonical}).
  Hence by \ref{Tomari's cor}, $\bp$ is the canonical weight for
  $F_t$ $(t\neq 0)$.
  Since $\bp(F_t)=\bp(\bl)$  for all $t\in C$,
  we can see that the deformation  $\pi:{\cal X}\to C$  defined by
  $\{F_t\}$
  admits the simultaneous canonical modification ${\cal X}(\bp)$
  in the same way as in the proof of \ref{deform 1}.
  Therefore
  $\gamma_m(X,0)=\gamma_m({\cal X}_1,0)$.
\end{pf}

\begin{exmp}
  (Watanabe)
  One can see in \cite{W2}
95-examples of deformations such as in Proposition \ref{deform 2}.
  For example,
 let  $X\subset \bC ^4 $  be defined by
$f=x_0^2+x_1^3+x_2^7+x_3^{43+s}+x_0x_1x_2x_3=0$
 $(s\geq 0)$
 and  $Y\subset \bC ^4 $ by $g=x_0^2+x_1^3+x_2^7+x_3^{42}=0$.
  Let  $\bp$  be $(21, 14, 6, 1)$,
  then  $X(\bp) \to X$  is the canonical modification and  $\bp/42$  is the
weight of
  the quasi-homogeneous polynomial  $g$.
  One can construct a family $\{F_t\}$ connecting $f$   and $g$  as in the
proof of
  Proposition \ref{deform 2}.
\end{exmp}

\section{A counter example to the conjecture}
\begin{say}
In this section
  we show a counter example to the conjecture written in the introduction.
  Let  $f$  be the polynomial:
  $x_0x_1x_2x_3+\alpha x_0^3 +\beta x_1^2x_2^2 +\gamma x_1^{a_1} +
  \delta x_2^{a_2}+\epsilon x_3^{a_3} \in \bC[x_0,\ldots ,x_3]$,
  with $a_i\geq 6$  and $\alpha,\ \beta,\ \gamma, \ \delta, \ \epsilon \in
  \bC$ general.
  Then  $f$  is non-degenerate and defines an isolated, non-canonical,
  log-canonical singularity $(X,0)$  at  the origin by
  \ref{characterization of log-canonical}.
  The essential cone is as follows:
$$\Cf=
    \{\bq\in \sigma |2q_0-q_1-q_2-q_3\geq 0,\ -q_0+q_1+q_2-q_3\geq 0,
   \ (a_i-1)q_i-\sum_{j\neq i}q_j\geq 0,\ i=1,2,3\}$$
Here $\Cf$ has no absolutely minimal vector.
  In fact,
  it is easy to see that  $(2,2,1,1) $  and $(2,1,2,1)$ belong to
  $\Cf$  but neither $(2,1,1,1)$  nor (1,1,1,1) does.
  This shows that under these coordinates
  there is no weighted blow-up which is the canonical
  modification of  $(X,0)$  by \ref{main theorem}.
  In the following,
  we prove the same statement under arbitrary coordinates.
\end{say}

\begin{lem}
  If  $ Y\to X$  is the canonical modification of  $(X,0)$,
  then $-K_Y^3>3/2$.
\end{lem}
\begin{pf}
  We use the notation in \ref{quadrant} and \ref{notation star}.
  Denote the fan consisting of all faces of the
  positive quadrant in $N_{\bR}$ by $\D$.
  Let  $\bq$  be $(2,1,2,1)$,
  and $\varphi:T_N(\D(\bq))\to \bC^4$,
  $\varphi'=\varphi |_{X(\bq)}:X(\bq) \to X$
  be the $\bq$-blow-ups of $\bC^4$ and $X$  respectively
  under the given coordinates.
  First we prove that $X(\bq)$  has log-terminal singularities.
  For any resolution  $\psi:\tilde X \to X(\bq)$  of the singularities
  on $X(\bq)$,
  we can write  $K_{\tilde X}=\psi^*{\varphi'}^*(K_X)+\sum_i a_iE_i$
with $a_i\geq -1$  for all exceptional divisors  $E_i$,
  because  $(X,0)$  is log-canonical.
  On the other hand,
  by \ref{valuation of function} and \ref{birational canonical},
  $K_{T_N(\D(\bq))}+X(\bq)=\varphi^*(K_{{\bC}^4}+X)
  -D_{\bq}$,
  since  $\bq(\bl)-1-\bq(f)=-1$.
  Therefore if we write:
  $K_{\tilde X}
  =\psi^*((K_{T_N(\D(\bq))}+X(\bq))|_{X(\bq)})+\sum_jm_jE_j$,
  then $m_j>-1$ for every exceptional divisor  $E_j$  of $\psi$.
  Hence,  by \ref{normality},
  $X(\bq)$  has at worst log-terminal singularities.
  Note that there is a non-canonical singularity,
  because $m_j=-1/2$  for  $E_j$ which corresponds to the vector
  $(2,2,1,1)$.
  Next construct a flat deformation  $\pi:({\cal X},0)\to(C, 0)$
  by  $(1-t)f+t(x_0^3+x_1^6+x_2^3+x_3^6)$
  as in \ref{deform 2}
  so that  $({\cal X}_0,0)\simeq (X,0)$, and $({\cal X}_1,0)$  is
  defined by non-degenerate
  weighted homogeneous polynomial $x_0^3+x_1^6+x_2^3+x_3^6$
 with the weight
 $\dfrac{1}{6}(2,1,2,1)$.
  This deformation is proved to be an (FG)-deformation
  (for the definition cf. \cite{Isiml}) as follows:
  Let  $\Phi:{\cal X}(\bq)\to {\cal X}$  be
  the restriction of  $\varphi\times id_{\bC}$  onto
  the proper transform  ${\cal X}(\bq)$  of  ${\cal X}$  in
  $T_N(\D(\bq))\times \bC$;
  since  $\bq(f)=\bq(\bl)=\bq(F_t)$  for  $t\in C$,
  ${\cal X}(\bq)_t$  is the $\bq$-blow-up ${\cal X}_t(\bq)$
  of ${\cal X}_t$
  for every  $t\in C$
  as in the proof of  \ref{deform 1};
  here  ${\cal X}(\bq)_t$ has at worst canonical singularities for
  $t\neq 0$ by \ref{Tomari's cor}
  and ${\cal X}(\bq)_0=X(\bq)$  has at worst log-terminal singularities as
  proved above;
  on the other hand, it is clear that  $K_{{\cal X}(\bq)}$  and
  $K_{{ X}(\bq)}$ are both  $\bQ$-Cartier divisors;
  hence by Proposition 7 of \cite{Stevens},
  ${\cal X}(\bq)$  admits at worst canonical singularities;
  one can easily see that
  $K_{{\cal X}(\bq)}=-(D_{\bq}\times C)|_{{\cal X}(\bq)}$
  which is $\varphi$-ample,
  which shows that ${\cal X}$ admits the canonical modification
  ${\cal X}(\bq)$ (i.e. $\pi$  is an (FG)-deformation).

  Now we can apply the upper semi-continuity theorem on $\{\gamma_m\}$
  (Theorem 1 of \cite{Isiml})
  to our (FG)-deformation  $\pi$.
  Since  $-K^3/3!$  of the canonical modification is the coefficient
  of the leading term of a function $\gamma_m$ in $m$,
  it follows
  that $-K_Y^3\geq -K_{{\cal X}(\bq)_t}^3=\sum q_i/\Pi q_i =3/2$
  for  $t\in C-\{0\}$.
  Here the equality does not hold.
  Because if it does,
  $\pi$ would admit the simultaneous canonical modification
  $\Psi:{\cal Y}\to {\cal X}$  by Corollary 1.11 on \cite{Isiml}.
  Since the simultaneous canonical modification must be the canonical
  modification of  ${\cal X}$  by  \cite{Stevens} again,
  ${\cal Y}$  would coincide with  ${\cal X}(\bq)$.
  However  ${\cal X}(\bq)_0$  has a non-canonical singularity
  as is seen above.
\end{pf}
\begin{say}
  Now we assume that there are coordinates $y_0,\ldots ,y_3$  on $\bC^4$
  and a weight  $\bp=(p_0,\ldots ,p_3)$  such that the $\bp$-blow-up
  $X({\bp})\to X$  under these coordinates gives the canonical
  modification, and then will induce a contradiction.
  Let  $g(y)=0$ be the defining equation of  $X$  under these coordinates.
  By \ref{main theorem}, it follows that $\bp(g)=\bp(\bl)$
   and therefore $-K^3_{X(\bp)}=\sum_ip_i/\Pi _i p_i>3/2$.
  Now it is easy to prove that at least three of
  the  $p_i$'s must be  1.
  Write the coordinates transformation as follows:
$$(T_i)\ \ \ \ \ \ \ \ \ \ \  x_i=\sum _{m\in
   \bZ_{\geq 0}^4}a_m^{(i)}y^m\ \
   \ \ (a_m^{(i)}\in \bC).$$
  We may assume that the coefficient of  $y_i$  in $(T_i)$  is not
  zero for
  each $i$ by reordering $\{y_i\}$'s.
  Then $y_0^3\in g$ (see  \ref{coefficient of divisor} for the notation),
  since  $x_0^3\in f$ and this is the unique monomial of degree 3 in $f$.
 Therefore $\bp(3,0,0,0)\geq \bp(\bl)$ which means $p_0\geq 2$,
  since  $\bp$  must be in $\Cg$ by \ref{main theorem}.
  Then one obtains the fact that
$a_{0,1,0,0}^{(0)}=a_{0,0,1,0}^{(0)}=a_{0,0,0,1}^{(0)}
  =0$,
  otherwise $y_i^3\in g$, for  $i=1,2,3$  which induce $\bp(0,3,0,0)\geq
  \bp(\bl)$ and so on,  therefore $3\geq p_0+3$  a contradiction.
  One can also prove that $a_{0,0,1,0}^{(1)}= a_{0,1,0,0}^{(2)}=0$
  in the same way.
  Then it follows that  $y_1^2y_2^2\in g$,
  because this monomial comes from the term  $x_1^2x_2^2 $  and is not
  cancelled
  by the contribution from other terms.
Hence  $\bp$  must satisfy  $\bp(0,2,2,0)\geq \bp(\bl) $
  which is equivalent to $4\geq p_0+3$, a contradiction.
\end{say}


\makeatletter \renewcommand{\@biblabel}[1]{\hfill#1.}\makeatother

\vskip 1truecm
\hskip 8truecm
Department of Mathematics

\hskip 8truecm
Tokyo Institute of Technology

\hskip 8truecm
Oh-okayama Meguro,

\hskip 8truecm
152 Tokyo, Japan

\hskip 8truecm
e-mail address: shihoko@@math.titech.ac.jp


\begin{thebibliography}{11}

\bibitem{F} A.R. ~Fletcher, {\em Plurigenera of 3-folds and weighted
  hypersurfaces},  thesis submitted for the degree of Ph.D. Univ. Warwick
  (1988)




\bibitem{IGor} S.~Ishii, {\em On isolated Gorenstein singularities},
   Math. Ann. {\bf 270} (1985), 541-554

\bibitem{def} \bysame, {\em Small deformations of normal singularities}
   Math. Ann. {\bf 275} (1986), 139-148


\bibitem{Iasym} \bysame, {\em The asymptotic behavior of pluri-genera for a
normal isolated
singularity}, Math. Ann. {\bf 286} (1990), 803-812

\bibitem{Isiml} \bysame, {\em  Simultaneous canonical modifications of
deformations of isolated
singularities},  Algebraic Geometry and Analytic Geometry,  Proceeding of the
Satellite Conference of ICM 90, Springer Lecture Note, (1991), 81-100

\bibitem{IF} \bysame, {\em On Fano 3-folds with non-rational
   singularities and two-dimensional base}, Abh. Math. Sem. Univ. Hamburg
  {\bf 64} (1994), 249-277


\bibitem{IB} \bysame, {\em  The weighted blowing ups of singularities
    with $\kappa = 0$},
  to appear in the Proceeding of the Conference of Singularities
  Beijin  (1994)


\bibitem{Kouchinirenco} A.G.~Kouchinirenco, {\em
  Poly{\'e}dres de Newton et nombre de Milnor},
  Inventiones math. {\bf 32}, (1976) 1-31

\bibitem{L} J. ~Lipman, {\em Double point resolutions of deformations
  of rational
  singularities}, Comp. Math. {\bf 38}, (1979) 37-43

\bibitem{M} S. ~Mori,
  {\em  Flip theorem and the existence of minimal models for
  3-folds}, J. Amer. math Soc. {\bf 1}, (1988) 117-253

\bibitem{Oda75} T. ~Oda, {\em Lectures on torus embeddings and applications},
  (Based on joint work with Katsuya Miyake), Tata Inst. of Fund. Research
{\bf 58},
   Springer-Verlag, Berlin-Heidelberg-New York (1978)

\bibitem{Reid} M. ~Reid, {\em Canonical 3-folds}, Journ\'ees de
  G\'eom\'etrie Algebrique d'Angers, A. Beauville, editor,
  Sijthoff and Nordhoff, Alphen aan den Rijn (1980) 273-310

\bibitem{Reid2} \bysame, {\em Young person's guide to
  canonical singularities}, Proceedings of Symposia in Pure Math. {\bf 46},
  (1987) 345-414


\bibitem{Stevens} J. ~Stevens, {\em  On canonical singularities as total
 spaces of deformations}, Abh. Math. Sem. Univ. Hamburg {\bf 58},
 (1988) 275-283

\bibitem{T} M. ~Tomari,
  {\em The canonical filtration of higher dimensional purely
  elliptic singularity of a special type},  Inventiones math. {\bf 104}
  (1991) 497-520


\bibitem{TW} \bysame, \& K. ~Watanabe, {\em
 On $L^2$-plurigenera of not-log-canonical Gorenstein isolated
 singularities},
 Amer. Math. Soc. {\bf 109}, (1990) 931-935

\bibitem{V} A.N. ~Varchenco, {\em Zeta-function of monodromy and
  Newton's diagram}, Inv. Math. {\bf 37}, (1976) 253-262


\bibitem{W} K. ~Watanabe, {\em
 On plurigenera of normal isolated singularities I},
Math. Ann. {\bf 250}, (1980) 65-94

\bibitem{W3} \bysame, {\em On plurigenera of normal isolated singularities
  II}, Advanced Studies in Pure math. {\bf 8} (1986) 671-685

\bibitem{W2} \bysame, {\em Table of standard polynomials of type (0,0)
 for 95-weights}, private letter (1994)


\bibitem{HW} \bysame, \& H.~Higuchi, {\em
 On certain class of purely elliptic singularities in dimension
$>2$}, Sci. Rep. Yokohama Nat. Univ. Sect. I, {\bf 30} (1983) 31-35



\bibitem{Y} T. ~Yonemura, {\em Hypersurface simple K3-singularities},
    Tohoku Math. J. Second Series {\bf 42} (1990) 351-380




\vskip 2truecm
\end{thebibliography}
\end{document}